\documentclass{revtex4}

\usepackage{indentfirst}
\usepackage{amsfonts}
\usepackage{amssymb}
\usepackage{graphicx}

\usepackage[reqno,intlimits]{amsmath}
\usepackage[polish,english]{babel}

\begin{document}

\title{A differential algorithm for the Lyapunov spectrum}
\author{Tomasz Stachowiak}
\email{toms@oa.uj.edu.pl}
\affiliation{Copernicus Center for Interdisciplinary Studies,\\
ul. Gronostajowa 3, 30-387 Krak\'ow, Poland}

\author{Marek Szydlowski}
\email{uoszydlo@cyf-kr.edu.pl}
\affiliation{Copernicus Center for Interdisciplinary Studies,\\
ul. Gronostajowa 3, 30-387 Krak\'ow, Poland\\}
\affiliation{Astronomical Observatory, Jagiellonian University\\
ul. Orla 171 30-244 Krak\'ow, Poland}

\begin{abstract}
We present a new algorithm for computing the Lyapunov exponents spectrum based
on a matrix differential equation. The approach belongs to the so called
continuous type, where the rate of expansion of perturbations is obtained for
all times, and the exponents are reached as the limit at infinity. It does not
involve exponentially divergent quantities so there is no need of rescaling or
realigning of the solution. We show the algorithm's advantages and drawbacks using
mainly the example of a particle moving between two contracting walls.
\end{abstract}

\maketitle

\section{Introduction}

Lyapunov characteristic exponents (LCE) measure the rate of exponential
divergence between neighbouring trajectories in the phase space. The standard
method of calculation of LCE for dynamical systems is based on the variational
equations of the system. However, solving these equations is very difficult or
impossible so the determination of LCE also needs to be carried out numerically
rather than analytically. 

The most popular methods which are used as an effective numerical tool to
calculate the Lyapunov spectrum for smooth systems relies on periodic
Gram-Schmidt orthonormalisation of Lyapunov vectors (solutions of the
variational equation) to avoid misalignment of all the vectors along the
direction of maximal expansion (\cite{Galgani}, \cite{Christansen}).

In some approaches, usually involving a new differential equation instead of
the variational one, the procedure of re-orthonormalisation is not used
\cite{Ranga} -- these are usually called the continuous methods. They are
usually found to be slower than the standard ones due to the underlying
equation being more complex than the variational one.
A comparison of various methods with and without orthogonalisation can be found
in \cite{Rama, Geist} and a recent general review in \cite{Skokos}.

The main goal of this paper is to present a new algorithm for obtaining the LCE
spectrum without the rescaling and realigning. This application is a
consequence of the equation satisfied by the Lyapunov matrix or operator (see
below) which was discovered in one of the authors' PhD thesis \cite{PhD}. The
particular numerical technique introduced here is the first attempt and is
open to further development so it still bears the disadvantages of the
usual continuous methods. However, in our opinion, the main advantages of the
approach lie in its founding equations and are as follows:
\begin{itemize}
\item The whole description of the LCE is embedded in differential geometry
from the very beginning, so that and it is straightforward to assign any metric
to the phase space including one with non-trivial curvature.
\item As the rate of growth is described by an operator (endomorphism) on the
tangents space, and the equation it satisfies is readily expressed with the
absolute derivative, the approach is explicitly covariant. (The exponents are
obviously invariants then, although their transformation properties 
still seem to be a live issue, see e.g. \cite{Eich}.)
\item There is no need for rescaling and realignment, as the main matrix
is at most linear with time, and encodes the full spectrum of LCE. 
\item Since we make no assumptions on the eigenvalues, there are no problems
with the degenerate case encountered in some other methods.
\item We rely on a single coordinate-free matrix equation, which reduces the
method's overall complexity.
\item The fundamental equation is not an approximation but rather the
differential equation satisfied by the so-called time-dependent Lyapunov
exponents. This opens potential way to analytic studies of the exponents.
\end{itemize}

It should be noted that the last points imply a hidden cost (in the current
implementation) of diagonalising instead of reorthonormalising, due to the
complex matrix functions involved. Fortunately, this procedure needs to be
carried out on symmetric matrices for which it is stable.

The natural domain of application of this method might be the General
Relativity and dynamical systems of cosmological origin -- already formulated
in differential geometric language \cite{Xin,Szyd1,Szyd2}. Of course this
still requires the resolution of the question of the time parameter, and
natural metric in the whole phase space (not just the configuration space
which corresponds to the physical space-time). Regardless of that choice,
however, the fundamental equation of our method will remain the same -- whether
one chooses to consider the proper interval as the time parameter, or find some
external time for an eight dimensional phase space associated to the four
dimensional space-time. This stems from the fact that our approach works on any
manifold.

Here, we wish to focus on the numerical aspect of the method, providing the
rough first estimates of its effectiveness. This is a natural question, after
the theoretical motivation for a given method has been established, namely how
well it performs numerically. There are obviously many ways of
translating the method into code, and we hope for future improvements,
nevertheless, the presented implementation can be considered a complete,
ready-to-use tool. In the next sections we review the derivation of the main
equation and then proceed to the simple mechanical examples for testing and
results.

\section{The base system of equations}

For a given system of $n$ ordinary differential equations
\begin{equation}
    \dot{x} = f(x,t),\;\; x\in\mathbb{R}^n,
    \label{dyn_sys}
\end{equation}
the variational equations along a particular solution $\varphi(t)$ are defined
as
\begin{equation}
    \dot{Z} = A Z,\; A:= \left.\frac{\partial f}{\partial x}\right|_{x=\varphi(t)},
    \label{var_eq}
\end{equation}
and the largest Lyapunov exponent can be defined as
\begin{equation}
    \lambda_{\mathrm max} = \lim_{t\rightarrow\infty} \frac{\ln(\|Z\|)}{t},
\end{equation}
for almost all initial conditions of $Z$. From now on we take the norm to be
\begin{equation}
    \|Z\| := \sqrt{Z^{T}Z},
\end{equation}
where $Z$ is treated as a column vector, and $T$ denotes transposition. That is
to say the metric in the tangent space is Euclidean, as is usually assumed
for a given physical systems. This needs not be the case, and a fully covariant
derivation of the main equation can be found in \cite{PhD}.

The above definitions are intuitively based on the fact
that for a constant $A$, the solution of \eqref{var_eq} is of the form
\begin{equation}
    Z = \exp(A t)Z(0),
\end{equation}
and $\lambda_{\mathrm{max}}$ is the greatest real part of the eigenvalues of
$A$. In the simplest case of a symmetric $A$, the largest exponent is exactly
the largest eigenvalue. To extend this to the whole spectrum, we note that any
solution of
\eqref{var_eq} is given in terms of the fundamental matrix $F$ so that
\begin{equation}
    Z = F Z(0),\;\; \dot{F} = AF,\;\; F(0)=\mathbf{1}.
\end{equation}
Then (if the limit exists), the exponents are
\begin{equation}
    \{\lambda_1,\ldots,\lambda_n\} =
    \mathrm{spec}\left(\lim_{t\rightarrow\infty}\frac{\ln(F^TF)}{2t}\right).
\end{equation}
Since $F^TF$ is a symmetric matrix with non-negative eigenvalues, the logarithm
is well defined. The additional factor of 2 in the denominator results from the
square root in the definition of the norm above.

As we expect $F$ to diverge exponentially, there is no point in integrating the
variational equation in itself, but rather to look at the logarithm. To this
end we introduce the two matrices $M$ and $L$:
\begin{equation}
    M := FF^T = F(F^TF)F^{-1} =: \exp(2L).
\end{equation}
Clearly $M$ has the same eigenvalues as $F^TF$ to which it is connected by a
similarity transformation, and the eigenvalues of $L$ behave as $t\lambda$ for
large times
\begin{equation}
    \{\lambda_1,\ldots,\lambda_n\} = 
    \mathrm{spec}\left(\lim_{t\rightarrow\infty}\frac{L}{t}\right).
\end{equation}
That is why we call $L$ the Lyapunov matrix.

To derive the differential equation satisfied by $L$ we start with the
derivative of $M$ and use the property of the matrix (operator) exponential
\begin{equation}
\begin{aligned}
     \dot{M} = AFF^T + FF^TA^T &= 
     e^{2L}\int_0^1 e^{-2aL}2\dot{L} e^{2aL}\mathrm{d}a\\
     e^{-2L}Ae^{2L}+A^T &= 2\int_0^1 e^{-2a[L}\dot{L}\mathrm{d}a,
\end{aligned}
\end{equation}
where we have introduced a concise notation for the the adjoint of $L$ acting
on any matrix $X$ as
\begin{equation}
    [L(X) = [L,X],
\end{equation}
and used its property
\begin{equation}
    e^{-L}Xe^{L} = e^{-[L}X.
\end{equation}
Next, the integral is evaluated taking the integrand as a formal power series
in $[L$
\begin{equation}
\begin{aligned}
    e^{-2[L}A + A^T = 2\frac{e^{-2[L}-\mathbf{1}}{-2[L}\dot{L},
    \label{integrated}
\end{aligned}
\end{equation}
where the fraction is understood as a power series also, so that there are in
fact no negative powers of $[L$. Alternatively one could justify the above by
stating that the function
\begin{equation}
    \psi_1(u) = \frac{1-e^{-2u}}{u}
\end{equation}
is well behaved on the spectrum of $[L$ which is contained in $\mathbb{R}$. As
$\psi_1$ is never zero for a real argument, we can invert the operator on the
left-hand side of \eqref{integrated} to get
\begin{equation}
    \frac{[L}{1-e^{-2[L}}\left(e^{-2[L}(\theta+\omega)+\theta-\omega\right)
    = \dot{L},
\end{equation}
where the symmetric and antisymmetric parts of $A$ are
\begin{equation}
    \theta = \frac12(A+A^T),\;\; \omega = \frac12(A^T-A^T).
\end{equation}
This allows for the final simplification to
\begin{equation}
    \dot{L} = [L \coth\left([L\right)\theta - [L,\omega]
    \label{main_eq}
\end{equation}
The function $\psi_2(u) = u\coth(u)$ should be understood as the appropriate
limit at $u=0$, so that it is well behaved for all real arguments. As was
proven in \cite{PhD}, the above equation is essentially the same in general
coordinates:
\begin{equation}
    \nabla_v L = \psi_2([L)\theta - [L,\omega],
\end{equation}
where $v$ is the vector field associated with  \eqref{dyn_sys}, and $\nabla$ is
the covariant derivative.

Note that in this form it is especially easy to obtain the known result for the
sum of the exponents. Since trace of any commutator is zero, the only term left
is the ``constant'' term of $\psi_2$ which is 1 (or rather the identity
operator) so that
\begin{equation}
    \frac{\mathrm{d}}{\mathrm{d}t}\sum_{i=1}^nt\lambda_i =
    \frac{\dot{V}}{V} = \mathrm{tr}\,\theta,
\end{equation}
where $V$ is the volume of the parallelopiped formed by $n$ independent
variation vectors.

Another simple consequence occurs when the $\theta$ matrix is zero, the whole
equation becomes a Lax equation
\begin{equation}
    \dot L = [\omega,L],
\end{equation}
which preserves the spectrum over time, so that $L/t$ tends to zero at
infinity. Another way of looking at it is that it is a linear equation in $L$
and the matrix of coefficients $[\omega$ is antisymmetric in the adjoint
representation, so that the evolution is orthogonal and the matrix norm
(Frobenius norm to be exact) of $L$ is constant which means $L/t$ tends to the
zero matrix. The simplest example of this is the harmonic oscillator or any
critical point of the centre type. The variations are then vectors of constant
length and the evolution becomes a pure rotation. The authors are not aware of
any complex or non-linear system that would exhibit such simple behaviour.
Already for the mathematical pendulum such picture is achieved only
asymptotically for solutions around its stable critical point. One could
expect that a system with identically zero exponents might not be
``interesting'' enough to incur this kind of research.

We have thus arrived at a dynamical system determining $L$, with right-hand
side being given as operations of the adjoint of $L$ on time dependent (through
the particular solution) matrices $\theta$ and $\omega$. The next section deals
with the practical application of the above equation.

\section{Exemplary implementation}

The main difficulty in using \eqref{main_eq} is the evaluation of the function
of the adjoint operator. Since we will be integrating the equation to obtain
the elements of the matrix $L$, it would be best to have the right-hand side as
an explicit expression in those elements. This can be done for the $n=2$ case,
but already for $n=3$ one has dozens of terms on the right, and for higher
dimensions the number of terms is simply too large for such an approach to be
of practical value. An alternative (although equivalent) dynamical system
formulation for the mentioned low dimensionalities have been studied in
\cite{Ranga}, but again the complexity of the equations increases so fast with
the dimension that the practical value is questionable. Our method, on the
other hand, can be made to rely on the same equation for all dimensions, and the
only complexity encountered will be the diagonalisation of a symmetric matrix
of increasing size.

Another problem lies in the properties of the $\psi_2(u)$ function which,
although finite for real arguments, has poles at $u=\pm i\pi$. This means that
a series approximation is useless, as it would converge only for eigenvalues
smaller
than $\pi$ in absolute value, whereas we expect them to grow linearly with time
and need the results for $t\rightarrow\infty$. On the other hand,
$\psi_2(u)\rightarrow|u|$ for large $u$ but, unfortunately, the adjoint
operator always has $n$ eigenvalues equal to zero, and for Hamiltonian systems
it is also expected that two eigenvalues tend to zero.

Thus, as we require the knowledge of $\psi_2$ for virtually any symmetric
matrix $L$, and we are going to integrate the equation numerically anyway,
we will resort to numerical method for this problem. Because the matrix
$L$ is symmetric (Hermitian in an appropriate setting) so is its adjoint
$[L$, and the best numerical procedure to evaluate its functions is by direct
diagonalisation \cite{Higham}. Obviously this is the main disadvantage of the
implementation method as even for symmetric matrices, finding the eigenvalues
and all
the eigenvectors is time-consuming. So far the authors have only been able to
find one alternative routine which is to numerically integrate not $L$
itself but rather the diagonal matrix of its eigenvalues and the accompanying
transformation matrix of eigenvectors. However, due to the increased number of
matrix multiplications the latter method does not seem any faster than the
former.

With this in mind, let us see how the diagonal form of $L$ simplifies the
equation. First, we need to regard $[L$ as an operator, and since it is acting
on matrices we will adopt a representation where any $n\times n$
matrix becomes a $n^2\times 1$ matrix, i.e. a $n^2$ element vector constructed
by writing all the elements of successive rows as one column. $[L$ is then
a $n^2\times n^2$ matrix. Fortunately, one does not need to diagonalise $[L$
but only $L$ itself. As can be found by direct calculation, the eigenvalues of
the adjoin are all the differences of the eigenvalues of $L$. For example
\begin{equation}
    L_D = \left(\begin{matrix} L_1 & 0 \\ 0 & L_2 \end{matrix}\right)
    \;\Rightarrow\;
    [L_D = \left(\begin{matrix} 
    0 & 0 & 0 & 0\\
    0 & L_1-L_2 & 0 & 0\\
    0 & 0 & L_2-L_1 & 0\\
    0 & 0 & 0 & 0 \end{matrix}\right),
\end{equation}
where the subscript $D$ denotes ``diagonal''.

Now let us assume we also have the transformation matrix such that
\begin{equation}
    L = Q L_D Q^T.
\end{equation}
Then, instead of bringing the whole equation to the eigenbasis of $L$, one can
only deal with the $\theta$ matrix in the following way
\begin{equation}
    \psi_2([L) \theta = Q\, \psi_2([L_D)Q^T\theta =
    Q\, \psi_2([L_D)\widetilde{\theta}\,Q^T,
    \;\; \widetilde{\theta}:=Q^T\theta\,Q.
    \label{eig_tr}
\end{equation}
Of course, the other term of equation \eqref{main_eq} can be evaluated as the
standard commutator. For the above example the $\psi_2$ part would be
\begin{equation}
    \psi_2([L_D) = \left(\begin{matrix}
    1 & 0 & 0 & 0\\
    0 & \psi_2(L_1-L_2) & 0 & 0\\
    0 & 0 & 1 & 0\\
    0 & 0 & 0 & \psi_2(L_2-L_1) \end{matrix}\right),
\end{equation}
and converting $\widetilde{\theta}$ to a vector we get
\begin{equation}
    \psi_2([L_D) \tilde\theta = 
    \psi_2([L_D) \begin{pmatrix}\tilde\theta_{11}&\tilde\theta_{12}\\
    \tilde\theta_{21}&\tilde\theta_{22}\end{pmatrix} =
    \psi_2([L_D)
    \begin{pmatrix}\tilde\theta_{11}\\\tilde\theta_{12}\\
    \tilde\theta_{21}\\\tilde\theta_{22}\end{pmatrix} = 
    \begin{pmatrix}\tilde\theta_{11}\\\psi_2(\Delta_{12})\tilde\theta_{12}
    \\\psi_2(\Delta_{21})\tilde\theta_{21}\\\tilde\theta_{22}\end{pmatrix},
\end{equation}
which corresponds to the usual 2 by 2 matrix of
\begin{equation}
    \psi_2([L_D)\tilde\theta = \begin{pmatrix}
    \tilde\theta_{11}&\psi_2(\Delta_{12})\tilde\theta_{12}\\
    \psi_2(\Delta_{21})\tilde\theta_{21}&\tilde\theta_{22}\end{pmatrix},
\end{equation}
where $\Delta_{ij} := L_i-L_j$ are the differences of the eigenvalues of $L$.
In general the appropriate ($n\times n$) matrix elements read
\begin{equation}
    \left[\psi_2([L_D)\tilde\theta \right]_{ij} =
    \psi_2(\Delta_{ij})\tilde\theta_{ij},
\end{equation}
and this matrix needs to be transformed back to the original basis according to
\eqref{eig_tr} before being used in the main equation.

The matrix elements of $L$ will, in general, grow linearly with time. This is
of course a huge reduction when compared with the exponential growth of the
perturbations, but one might want to make them behave even better by taking time
into account with
\begin{equation}
    L = (t+1)\Lambda.
    \label{rescal}
\end{equation}
The reason for the additional 1 is that we will specify the initial conditions
at $t=0$ and we want to avoid dividing by zero in the numeric procedure and
the limit at infinity (if it exists) is not affected by this change. Of
course, in the case of the autonomous systems any value of $t$ can be chosen
as initial (at the level of the particular solution) but the non-autonomous
case might require a particular value, which can be dealt with in a similar
manner.

We now have
\begin{equation}
    \dot\Lambda = 
    \frac{1}{t+1}\left(\psi_2\left((t+1)[\Lambda\right)\theta - \Lambda\right)
    -[\Lambda,\omega].
\end{equation}
As for the initial conditions, the fundamental matrix $F$ is equal to the
identity matrix at $t=0$, so that $L(0) = 0$ which in turn implies
$\Lambda(0)=0$.

We are now ready to state the general steps of the proposed implementation.
Choosing a specific numerical routine to obtain the solution for a time step
$h$ at each $t$ these are

\begin{enumerate}
\item Obtain the particular solution $\varphi(t)$, calculate the Jacobian
matrix $A$ at $\varphi$ and, from it, the two matrices $\theta(t)$ and
$\omega(t)$. Start with $\Lambda(0) = 0$

\item Find the eigenvalues $\lambda_i$ and eigenvectors $Q$ of $\Lambda$
\item Transform $\theta(t)$ to $\tilde\theta=Q^T\theta\,Q$.

\item Compute the auxiliary matrix $\beta_{ij}=
[\psi((t+1)[\Lambda_D)\tilde\theta]_{ij} =
\psi_2((t+1)\Delta_{ij})\tilde\theta_{ij}$ 

\item Take the derivative to be $\dot\Lambda(t) = 
(Q\,\beta\,Q^T -\Lambda)/(t+1) -[\Lambda,\omega]$, 
and use it to integrate the solution to the next step $\Lambda(t+h)$.

\item Repeat steps 2--5 until large enough time is reached.

\item The ``time-dependent'' Lyapunov exponents at time $t$ are simply
$(1+\frac{1}{t})\lambda_i$.

\end{enumerate}
The relation in the last point stems from the rescaling \eqref{rescal} and, as
can be seen, is only important for small values of $t$.

\section{Examples and Comparison}

As suggested in the preceding section, one expects this method of obtaining
the Lyapunov spectrum to be relatively slow. In order to see that, we decided to
compare it with the standard algorithm based on direct integration of the
variational equation and Gram-Schmidt rescaling at each step \cite{Galgani}.
Although usually the rescaling is required after times of order 1, we
particularly wish to study a system with increasing speed of oscillations
and frequent renormalisation will become a necessity. In other words, we want
to give the standard method a ``head start'' when it comes to precision. 

For numerical integration we chose the modified midpoint method (with 4
divisions of the whole timestep $H=4h$), which has the
advantage of evaluating the derivative fewer times than the standard
Runge-Kutta routine with the same accuracy. Since we intended
a simple comparison on equal footing, we did not try to optimise
either of the algorithms and wrote the whole code in Wolfram's Mathematica due
to the ease of manipulation of the involved quantities and operations (e.g.
matrices and outer products).

The most straightforward comparison is for the simplest, i.e. linear dynamical
system, for which the main equation \eqref{dyn_sys} is the same as the
variational one \eqref{var_eq}, with constant matrix $A$. The exponents are
then known to be the real parts of the eigenvalues of $A$. To include all kinds
of behaviour, we took a matrix which has a block form
\begin{equation}
    A = \begin{pmatrix} 
    1 & -8 & 0 & 0 \\
    12 & 1 & 0 & 0\\
    0 & 0 & \frac12 & 2\\
    0 & 0& 1& \frac12
    \end{pmatrix},
\end{equation}
whose eigenvalues are
$\{1+4i\sqrt{6},1-4i\sqrt{6},\frac12+\sqrt{2},\frac12-\sqrt{2}\}$, so that
the LCE are approximately equal to $\{1,1,1.91421,-0.914214\}$.

\begin{figure}[h]
\includegraphics[width=.7\textwidth]{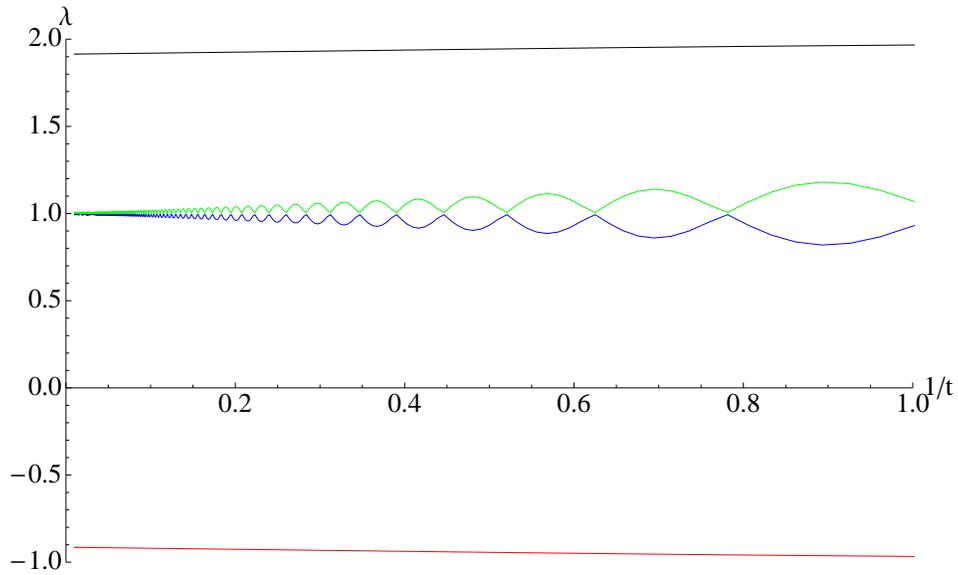}
\caption{The Lyapunov exponents for the linear system obtained with the
presented method.}
\label{const_1}
\end{figure}

\begin{figure}[h]
\includegraphics[width=.7\textwidth]{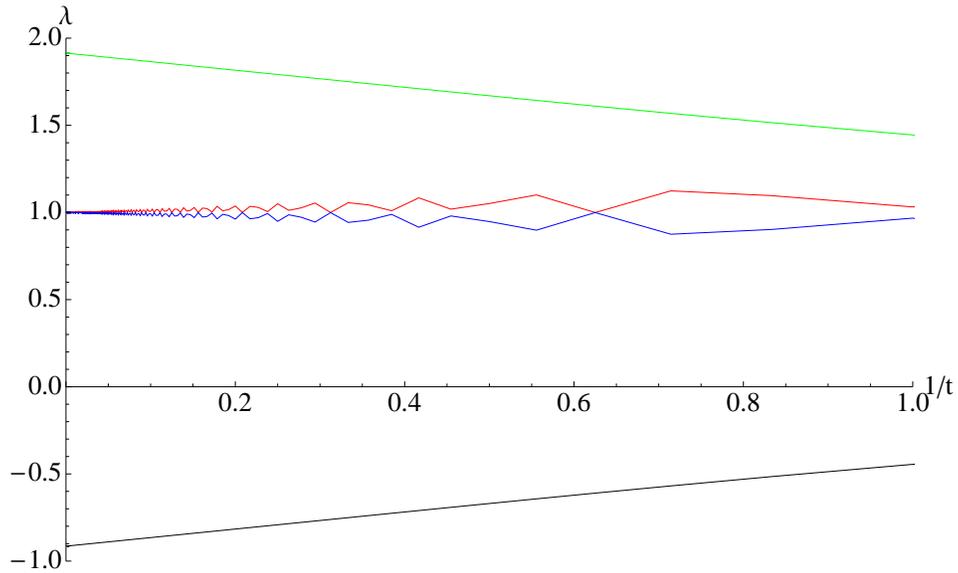}
\caption{The Lyapunov exponents for the linear system obtained with the
standard method.}
\label{const_2}
\end{figure}

The basic timestep was taken to be $h=0.01$ and the time to run from 0 to 1000.
The results are depicted in figures \ref{const_1} and \ref{const_2} with the
horizontal axis
representing the inverse time $1/t$ so that the sought for limit at infinity
becomes the value at zero which is often clearly seen from the trend of the
curves. We note that our method required 59 seconds, whereas the standard one
only took 17 seconds. The final values of LCE ($t=1000$) were
$\{1.00010,0.99990,1.91427, -0.91427\}$ and $\{0.99933, 0.99970, 1.91372,
-0.91372\}$, respectively. 
The shape of the curves is different, which is to be
expected because the matrix $\Lambda$ measures the true growth of the variation
vectors at each point of time, and the other method provides more and more
accurate approximations to the limit values of the spectrum.

Before we go on to the central example, which explores a system with accelerating
oscillations, we will present numerical results for the system which is
synonymous with chaos, namely the Lorenz system \cite{Lorenz}.

\begin{figure}[h]
\includegraphics[width=.7\textwidth]{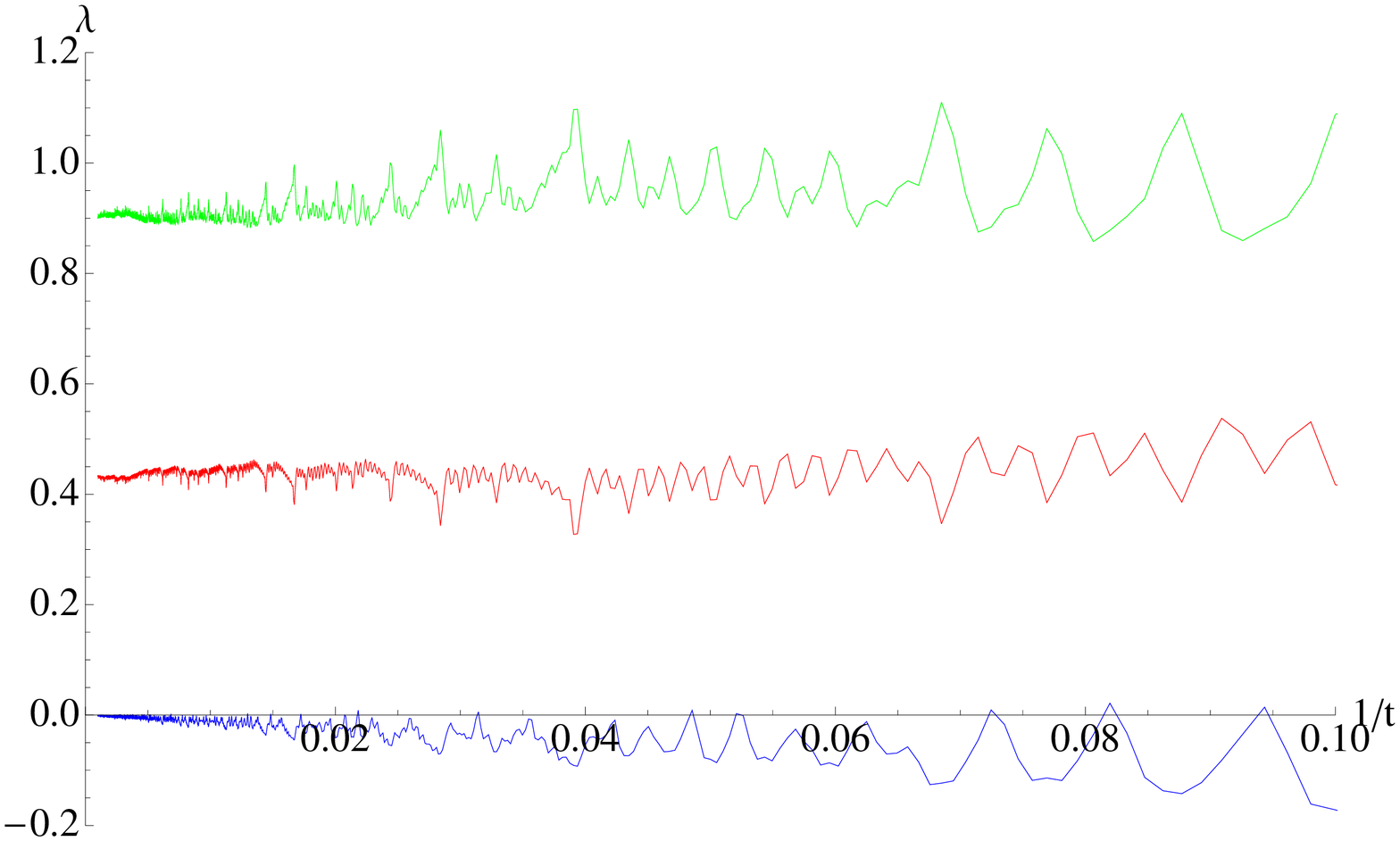}
\caption{The Lyapunov exponents for the Lorenz system obtained with the
new method. The middle exponent has been shifted by +15 to show the details
with a better scale.}
\label{lorenz_1}
\end{figure}

\begin{figure}[h]
\includegraphics[width=.7\textwidth]{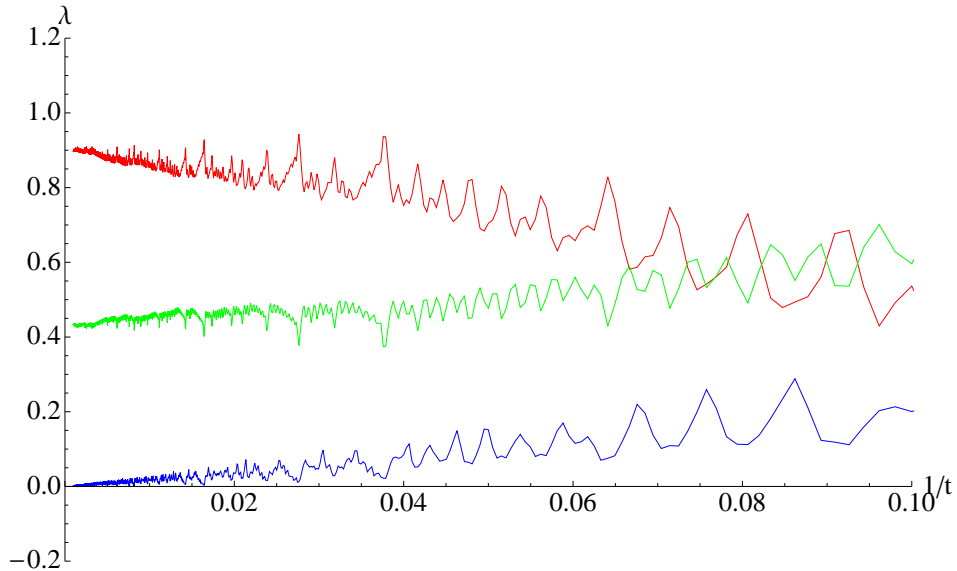}
\caption{The Lyapunov exponents for the Lorenz system obtained with the
standard method. The middle exponent has been shifted by +15 to show the
details with a better scale.}
\label{lorenz_2}
\end{figure}

The equations read
\begin{equation}
\begin{aligned}
    \dot x &= \sigma(y-x)\\
    \dot y &= x(\rho-z)-y\\
    \dot z &= xy-\beta z,
\end{aligned}
\end{equation}
where we took $\sigma = 10$, $\beta = 8/3$ and $\rho = 28$, and integrated the
equation for the initial conditions of $x=y=z=10$ from $t=0$ to 1000. Next we
integrated the respective methods for the exponents with the timestep of
$0.001$. Our method took about 591 seconds and the result is shown in figure
\ref{lorenz_1} with the final value of the spectrum $\{-14.5676, -0.0010078,
0.901931\}$. Note that we have
shifted the lowest exponent by $+15$, so that all three could be presented on
the same plot with enough detail. The standard method took about 152 seconds
and its outcome is shown in figure \ref{lorenz_2} with the final values of
$\{-14.5664, 0.00173051, 0.897989\}$ (we have shifted the graph in the same
manner as before).

One could note that there is less overall variation of the time-dependent
exponents in our method similarly to the previous example. A good estimate of
precision would be to calculate the sum of the exponents, which, in this
system, should be exactly equal to -41/3. The difference between that and the
numerical estimates were: for the standard method (at $t=1000$) 
$-9.04\times 10^{-7}$, and for our method $-9.67\times10^{-11}$ -- a
much better result. This seems to be the usual picture for the continuous
methods which trade computing time for precision.

A presentation of some more complicated, including both integrable and chaotic,
examples can be
found in \cite{PhD}, and for such systems also, we observe the concordance of
final results and the speed discrepancy. In order to see how one could benefit
from the new method we have to turn to another class of models, ones for which
the exponents present oscillatory behaviour. We found that for artificial
systems with accelerated oscillations our method performs better and present
here a simple physical model which exhibits such property.

\begin{figure}[h]
\includegraphics[width=.7\textwidth]{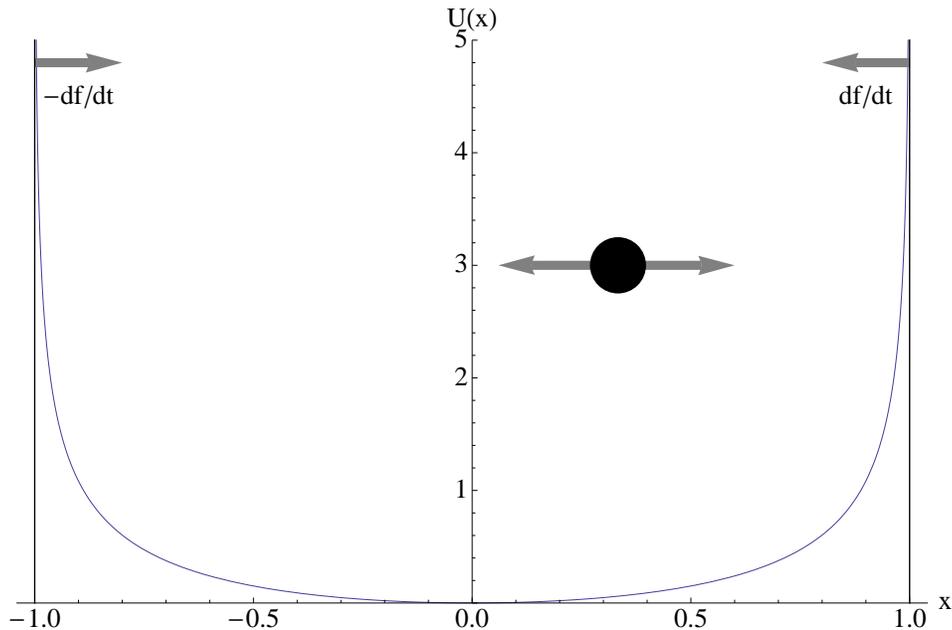}
\caption{The potential modelling two contracting walls. The vertical lines
correspond to the impenetrable barriers of ``infinite'' $U$. (The vertical
position of the ball has no meaning.)}
\label{potential}
\end{figure}

Consider a ball moving between two walls which are moving towards each other
and assume that the ball bounces off with perfect elasticity. As the distance
between the walls decreases and the speed of the ball increases it takes less
and less time for each bounce cycle. In order to model this analytically,
without resorting to infinite square potential well we will take the
following Hamiltonian system
\begin{equation}
    H = \frac12 p^2 + U\left(\frac{q}{f(t)}\right),\;\;
    U(x) = \frac12\mathrm{artanh}(x)^2,
\end{equation}
which depends on time explicitly via the function $f$ whose meaning is seen
as follows: the area hyperbolic tangent is infinite for $x=1$, so that the
potential becomes infinite for the position variable $q=f(t)$, so that $f$ is
simply half the distance between the walls. In particular we take it to be
\begin{equation}
    f(t) = \frac{t+20}{2t+20},
\end{equation}
so that it decreases from 1 to $\frac12$ slowing down but never stopping. The
reason for this is that we want the system not to end in a finite time, and
also that the worse behaved $f$ is the faster the numerical integration of the
main system itself will fail. The initial shape of the potential and the
systems setup is depicted in figure \ref{potential}

The slowing down of the walls, and the particular shape of $U$ allows us to
find rigorous bounds on the Lyapunov exponents. First we note, that the vector
tangent to a trajectory in the phase space, i.e. a vector whose components are
simply the components of $F$ from \eqref{dyn_sys}, is always (for any dynamical
system) a solution of the
variational equation. That means that just by measuring its length we can
estimate the largest exponent, since for almost all initial conditions the
resulting evolution is dominated by the largest exponent. Second,
the system is Hamiltonian and it must have two exponents of the same magnitude
but different signs, so analysing this particular vector will give us all the
information regardless of the initial condition. Thus we have to find
the following quantity
\begin{equation}
    \lambda(t) = \frac{1}{2t}\ln\left(\dot{q}^2 + \dot{p}^2\right)
    \label{v_norm}
\end{equation}
as a function of time, which boils down to finding bounds on the velocity and
acceleration of the bouncing ball.

As the Hamiltonian depends on time explicitly, the energy is not conserved, but
instead we have
\begin{equation}
    \dot{E} = \frac{\partial H}{\partial t} =
    U'\left(\frac{q}{f}\right)\left(-\frac{q\dot{f}}{f^2}\right).
\end{equation}
The velocity has its local maxima at $q=0$ when all the energy is in the
kinetic term, and we are lead to define a virtual maximal velocity by equating
the energy at any given time to a kinetic term
\begin{equation}
    \frac12 v_M^2 := E,
\end{equation}
we will take the positive sign of $v_M$, and assume it is non-decreasing as the
physical setup suggest. Differentiating the above we get
\begin{equation}
    v_M\dot{v}_M = - U'\left(\frac{q}{f}\right)q \frac{\dot{f}}{f^2}.
\end{equation}
Let us go back to the equation of motion for the momentum variable which reads
\begin{equation}
    \dot{p} = -U'\left(\frac{q}{f}\right)\frac{1}{f},
\end{equation}
and substitute that into the previous equation to get
\begin{equation}
    v_M\dot{v}_M = \dot{p}\, q \frac{\dot{f}}{f}.
\end{equation}
As mentioned above $f$ decreases very slowly at late times, which is when we
estimate the exponents anyway. We thus assume, that $\dot{f}$ is small enough
for the fraction on the right-hand side to be considered constant over one
cycle -- that is over the time $T$ in which the ball moves from the centre up
the
potential wall and back to the centre. This time will get shorter and shorter,
but also $\dot{f}$ will get closer and closer to zero. The standard problem of
the elastic ball and infinitely hard walls shows that the speed
transfer at each bounce is of the order of $\dot{f}$ so the (virtual) maximal
velocity will change as slowly as $f$ and we are entitled to average the
equations over one cycle:
\begin{equation}
    v_M \dot{v}_M = 
    \frac{1}{T}\frac{\dot{f}}{f}\int_t^{t+T}\dot{p}\, q \ \mathrm{d}t
    = \frac{1}{T}\frac{\dot{f}}{f}\left([qp]_t^{t+T} - 
    \int_t^{t+T} p^2\mathrm{d}t \right)
    \leq -\frac{\dot{f}}{f}v_M^2,
\end{equation}
where we integrated by parts and used the fact that at the beginning and end of
the cycle $q=0$, and also that the momentum is never greater than the maximal
velocity.

We are thus left with a bound in the form of a differential equation, and since
all the quantities involved are positive and non-decreasing ($-\dot{f} =
|\dot{f}|$), its solution will be the bound for the velocities
\begin{equation}
    \dot{q} \leq v_M \leq \frac{v_{M0}}{f},
\end{equation}
with the ``initial'' value of $v_M$ taken at sufficiently large $t$ so that
$\dot{f}$ is small.

Similar considerations can be carried out for the acceleration $\dot{p}$, only
this time we have to introduce a virtual points of return $q_R$, that is the
point at which all the energy is in the potential term
\begin{equation}
    U\left(\frac{q_R}{f}\right) := E,
\end{equation}
since at the real turning points the acceleration reaches its local maxima.
Note that this is not the same as $\dot{v}_M$ which describes the slow growth
of the consecutive maxima of $\dot{q}$ and not the maxima of its slope.

This definition allows us to express $q_R$ as a function of $v_M$ (via energy)
\begin{equation}
    q_r = f\tanh(v_M),
\end{equation}
and the acceleration is
\begin{equation}
\begin{aligned}
    \dot{p} \leq a_M &= \left|-U'\left(\frac{q_R}{f}\right)\right| = 
    \frac{f \mathrm{artanh}\left(\frac{q_R}{f}\right)}{f^2-q_R^2} =
    \frac{v_M}{f(1-\tanh^2(v_M))}\\
    a_M &= \frac{v_{M0}}{f^2}\cosh^2\left(\frac{v_{M0}}{f}\right).
\end{aligned}
\end{equation}

Bringing the two results together we see that
\begin{equation}
    \dot{q}^2+\dot{p}^2 \leq v_M^2 + a_M^2 < \infty,
\end{equation}
because the function $f$ does not tend to zero, and by \eqref{v_norm} both
Lyapunov
exponents must in turn be zero themselves. We also recognise that the
hyperbolic cosine factor in $a_M$ could produce a nonzero exponents if $f$
were to decrease to zero as $1/t$.

\begin{figure}[h]
\includegraphics[width=.7\textwidth]{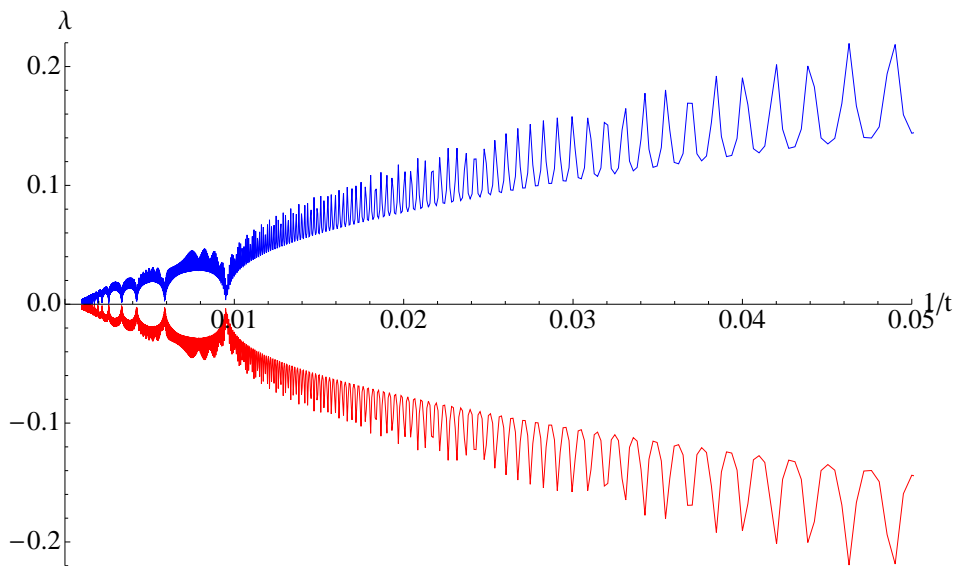}
\caption{The Lyapunov spectrum obtained with the new method for two contracting
walls (timestep $h=0.01$).}
\label{wall_1}
\end{figure}

\begin{figure}[h]
\includegraphics[width=.7\textwidth]{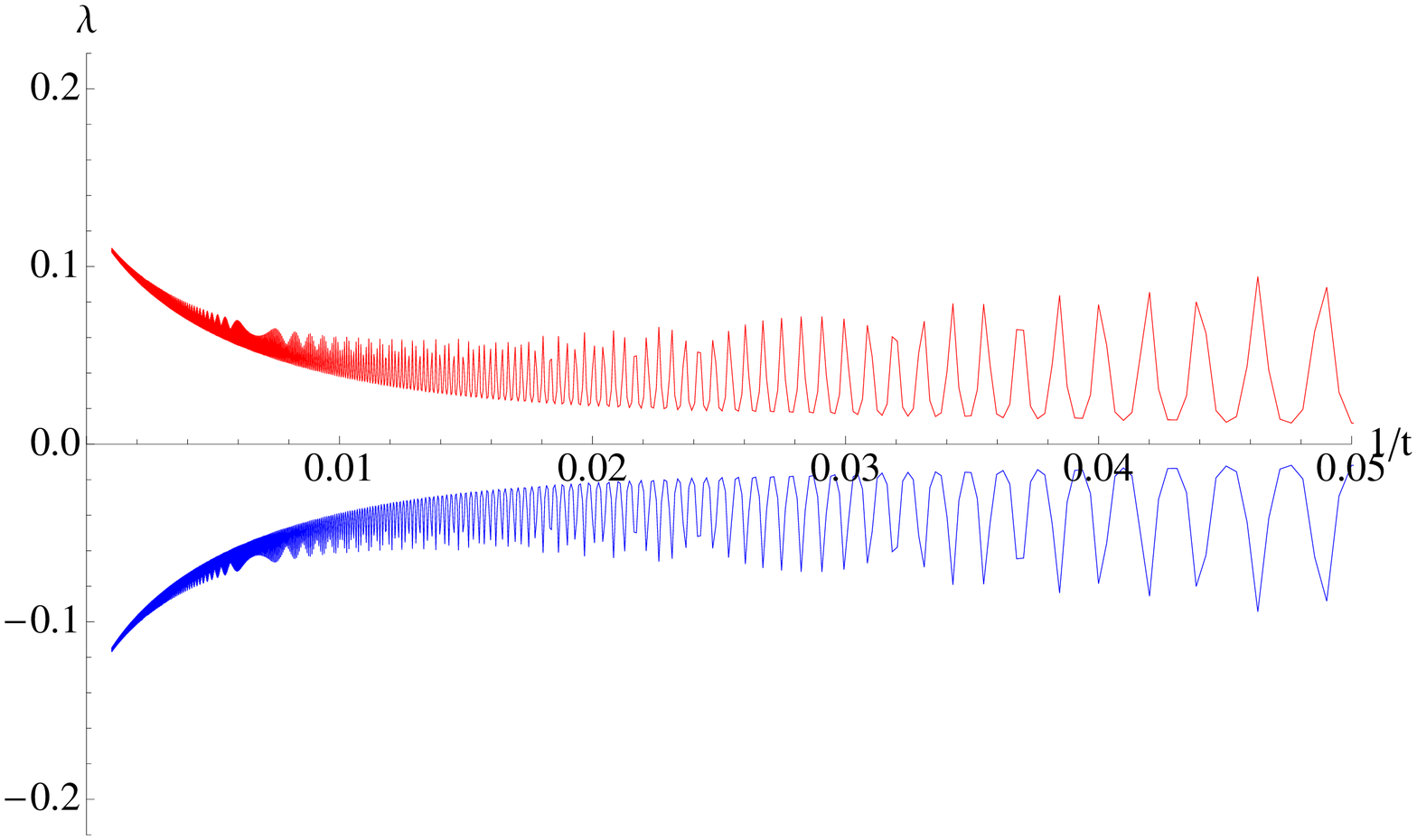}
\caption{The Lyapunov spectrum obtained with the standard method for two contracting
walls (time step $h=0.01$).}
\label{wall_2a}
\end{figure}

\begin{figure}[h]
\includegraphics[width=.7\textwidth]{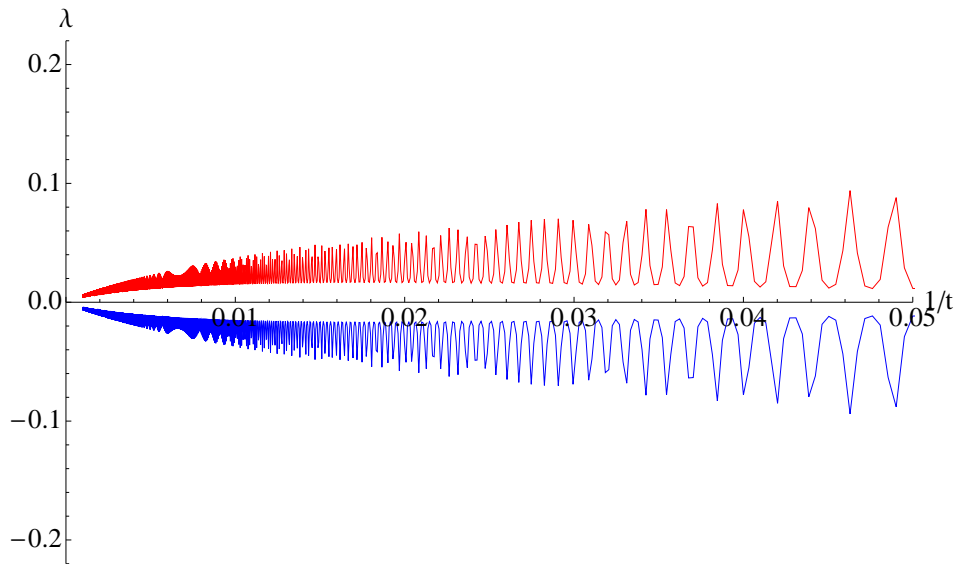}
\caption{The Lyapunov spectrum obtained with the standard method for two contracting
walls (time step $h=0.001$).}
\label{wall_2b}
\end{figure}

Let us now turn to the numerical results of both methods for this system. As
initial conditions we take $q(0)=0$ and $p(0)=1$. For
the same time step $h=0.01$ we see in figure \ref{wall_1} that the new method
predicts the values correctly, integrating for 59 seconds from $t=0$ to
$t=1000$. However for the standard one, as shown in figure
\ref{wall_2a}, we see the exponents diverging from zero, for the same time
limits, and integrating time of 17 seconds. It turns out $h=0.005$ also gives
divergent results and the correct behaviour is
recovered for $h=0.001$, shown in figure \ref{wall_2b}, for which the routine
takes 171 seconds.

\section{Conclusion}

We have presented a new algorithm for evaluation of the Lyapunov spectrum,
emerging in the context of differential geometric description of complex
dynamical systems. This description seems especially suitable for systems
found in General Relativity like, e.g., chaotic geodesic motion \cite{smerf}.
The main advantage of the base method is its covariant nature and concise,
albeit explicit, matrix equations that promise more analytic results in the
future. Also, this allows for study of curved phase spaces and general
dynamical systems -- not only autonomous or Hamiltonian ones.

The main differential equation, can be numerically integrated, giving a simple
immediate algorithm for the computation of the Lyapunov characteristic
exponents. It is in general slower than the standard algorithm (based on
Gram-Schmidt orthogonalisation), but the first numerical test suggest it works
betters in systems with increasing frequency of (pseudo-)oscillations. We show
this on the example of a simple mechanical system -- a ball bouncing between
two contracting walls.

Although in low dimensions the main equation can be cast into an explicit form
(with respect to the unknown variables), in general the numerical integration
requires diagonalisation at each step, which is the main disadvantage of the
method and the reason of its low speed. We hope to present a more developed
algorithm without this problem in the future.

\section{Acknowledgements}
This paper was supported by grant No. N N202 2126 33 of Ministry of Science and
Higher Education of Poland. The authors would also like to thank J. Jurkiewicz
and P. Perlikowski for valuable discussion and remarks.

\end{document}